\documentclass[10pt, a4paper]{amsart}

\usepackage[T1]{fontenc}
\usepackage[utf8]{inputenc}
\usepackage{amsmath, amssymb, amsfonts, amsthm}
\usepackage[margin=2.0cm]{geometry}
\usepackage[]{setspace}%\onehalfspacing
\usepackage[usenames,dvipsnames]{color}
\usepackage[, unicode]{hyperref}
\hypersetup{
    urlcolor=blue,
    colorlinks=true,
   linkcolor= blue,
   citecolor=blue
}
 \usepackage{comment, stmaryrd}
\usepackage[all]{xy}
 \usepackage{libertine}
 \usepackage[italic]{mathastext} 
\usepackage{bbm}
\usepackage{enumitem}
\setlist[enumerate]{label={\bf(\roman*)}, itemsep=1ex,leftmargin=1.4cm,topsep=1ex}
\setlist[enumerate,2]{label={\bf\alph*}., itemsep=1ex,leftmargin=0.5cm,topsep=1ex}
\setlist[enumerate,3]{label={\bf\roman*}., itemsep=1ex,leftmargin=0.5cm,topsep=1ex}

\theoremstyle{plain}
\newtheorem{theorem}{Theorem}[section]
\newtheorem{proposition}[theorem]{Proposition}
\newtheorem{lemma}[theorem]{Lemma}
\newtheorem{corollary}[theorem]{Corollary}

\newtheorem{fact}[theorem]{Fact}

\theoremstyle{definition}
\newtheorem{example}[theorem]{Example}
\newtheorem{definition}[theorem]{Definition}

\theoremstyle{remark}
\newtheorem{remark}[theorem]{Remark}

\newcommand\I{{\bf(i)}}
\newcommand\II{{\bf(ii)}}
\newcommand\III{{\bf(iii)}}

\newcommand\B{\mathrm{B}}
\newcommand\Lring{\mathfrak{L}_{\mathrm{ring}}}
\newcommand\Lval{\mathfrak{L}_{\mathrm{val}}}
\newcommand\Lac{\mathfrak{L}_{\mathrm{ac}}}
\newcommand\Lan{\mathfrak{L}_{\mathrm{an}}}

\newcommand\ps[1]{#1(\!(t)\!)}
\newcommand\ips[1]{#1[\![t]\!]}
\newcommand\UU{\mathcal{U}}
\newcommand\WW{\mathcal{W}}

\newcommand\GGG{\mathbf{G}}
\newcommand\res{\mathrm{res}}

\newcommand\Orb{\mathrm{Orb}}
\newcommand\perfect{Suppose $F$ is perfect.}

\renewcommand\mid{:}

\theoremstyle{plain}

\title[One-dimensional $F$-definable sets in $\ps{F}$]{\Large\rm One-dimensional $F$-definable sets in $\ps{F}$}
\author{Sylvy Anscombe}
\thanks{\today}
\address{Universit\'{e} Paris Cité and Sorbonne Universit\'{e}, CNRS, IMJ-PRG, F-75013 Paris, France}
\email{sylvy.anscombe@imj-prg.fr}

\begin{document}

\begin{abstract}
In this note we study one-dimensional definable sets in power series fields with perfect residue fields.
Using the description of automorphisms given by Schilling, in \cite{S44},
we show that such sets are unions of existentially definable in the language of rings,
allowing parameters.
We deduce that if $F$ is a perfect field of positive characteristic $p$, and $X$ is a subset of the $t$-adically valued  $\ps{F}$ that is definable in the language of valued fields with parameters from $F$,
then the subfield $(X)$ generated by $X$ is either contained in $F$ or equal to $F(\!(t^{p^{n}})\!)$, for some $n\geq0$.
The proof uses our earlier work on existentially definable subsets of henselian and large fields, of which power series fields are examples.
\end{abstract}

\maketitle
\setcounter{tocdepth}{2}
\begin{singlespacing}
\hypersetup{linkcolor=black}
% \tableofcontents
\hypersetup{linkcolor=blue}
\end{singlespacing}

\section{Introduction}

The theorem of Ax--Kochen/Ershov
(see e.g.~\cite[Theorem 3]{AK66})
provides an axiomatization of the first-order theory
of the power series field $\ps{F}$
in terms of the theory of the field of constants $F$,
in case that the characteristic of $F$ is zero.
On the other hand, if the characteristic of $F$ is positive, there is no known axiomatization.
Even in case that $F$ is finite, the model theory of $\ps{F}$ is largely unknown,
although see for example~\cite{Ku01,DS,AF16,Onay,Kostas_tame,ADF23}.
Neither is there a known description of the definable sets, with or without parameters.
Nevertheless, Ax--Kochen/Ershov-style results are known for some other valued fields
in positive characteristic,
for example
algebraically closed valued fields
\cite{Robinson},
separably closed value fields and those satisfying Kaplansky's hypothesis
\cite{Delon},
and more generally
tame and separably tame valued fields
\cite{Kuh16,KuhlmannPal}.
Recently there have also been results on deeply ramified and perfectoid fields
\cite{KR23,JK23}.

In this short note we study $\ps{F}$ in the case that $F$ is a perfect field of positive characteristic,
via the action of its group of automorphisms,
and by applying our previous analysis of existentially parameter-definable sets (\cite[Theorem 1]{A19}) in large fields that extended work of Fehm (\cite{Fehm}) on existentially definable sets in perfect large fields.
In Section~\ref{section:automorphisms}, we prove the following.

\begin{theorem}\label{thm:subfields_intro}
Let $F$ be a perfect field of characteristic $p>0$,
and let $X\subseteq \ps{F}$ be a subset definable in the language
$\Lval(F)$
of valued fields,
i.e.~allowing parameters from $F$.
Then either $X\subseteq F$ or there exists $n\in\mathbb{N}$ such that
$$
	(X)=F(\!(t^{p^{n}})\!),
$$
where $(X)$ denotes the subfield of $\ps{F}$ generated by $X$.
\end{theorem}

The basic principle underlying this theorem is that sets definable in a given structure are closed under automorphisms of that structure,
indeed even sets defined by types have this property.
Thus the theorem may be generalized to apply to any set $X$ defined in a ``suitable'' expansion of $\ps{F}$,
see Definition~\ref{def:suitable}
and Theorem~\ref{thm:subfields}.
\noindent We draw some further conclusions about subfields of $\ps{\mathbb{F}_{q}}$ and $\ps{\mathbb{F}_{q}}^{\mathrm{perf}}$ which are generated by $\Lval$-definable subsets, for any prime power $q$.

Finally, we can say something, albeit modest, about certain definable subsets of higher Cartesian powers of $\ps{F}$.
In Section~\ref{section:1-dim_orbits}, we prove the following.

\begin{theorem}\label{thm:1-dim_intro}
Let $F$ be a perfect field of characteristic $p>0$,
and let $\mathbf{a}=(a_{1},\dots,a_{n})$ be an $n$-tuple from $\ps{F}$ such that $F(\mathbf{a})/F$ is a field extension of transcendence degree $1$.
The orbit of $\mathbf{a}$ under the group of $\mathfrak{L}_{\mathrm{val}}(F)$-automorphisms
of $\ps{F}$
is
\begin{enumerate}
\item definable by an existential $\Lring$-formula with parameters from $F(t)$, and
\item definable by an $\Lring$-formula with parameters $F$.
\end{enumerate}
\end{theorem}

\section{Proof of Theorem~\ref{thm:subfields_intro}}
\label{section:first}

Throughout, $F$ will be a field
of characteristic $p>0$.
Our results are well-known in characteristic zero, following either directly from the classical theory of Ax--Kochen and Ershov,
or from other work such as \cite{JK10,Fehm}.
The field $\ps{F}$ of formal power series over $F$ in the indeterminate $t$
admits the {\em $t$-adic valuation}
$v_{t}:\ps{F}\twoheadrightarrow\mathbb{Z}\cup\{\infty\}$, which is discrete and maximal,
so in particular it is henselian.
The valuation ring of $v_{t}$ is $\mathcal{O}=\ips{F}$ and the maximal ideal is $\mathfrak{m}=t\ips{F}$.
We denote by $\UU$ the set of uniformisers of $v_{t}$,
i.e.~those elements of $\ps{F}$ of value $1$.
For $n\geq1$ we write $\WW_{n}$ for the set $t+t^{n}\ips{F}$,
Thus $\WW_{1}=\mathfrak{m}$ is the maximal ideal,
and we have the chain
\begin{align}\label{eq:chain_1}\tag{1}
		\mathcal{O}
	\supset
		\mathfrak{m}
	=
		\WW_{1}
	\supset
		\UU
	\supset
		\WW_{2}
	\supset
		\WW_{3}
	\supset
		\ldots
	\supset
		\bigcap_{n\geq1}\WW_{n}
	=
		\{t\}.
\end{align}
For any subset $A$ of a ring, and for each $l\in\mathbb{N}$, we denote by $A^{(l)}:=\{a^{l}\mid a\in A\}$ the set of $l$-th powers of elements of $A$.

\subsection{A simple application of Hensel's Lemma}
\label{section:polynomial}

We begin with a simple application of Hensel's Lemma.

\begin{lemma}\label{lem:Hensel_2}
Suppose $F$ is perfect.
Let
$a=\sum_{i=0}^{\infty}a_{i}t^{i}\in\ips{F}$
and let $m\geq2$.
There exists a
polynomial
$f\in F[X]$
and 
$s\in\WW_{m}$
such that
$a=f(s)$.
In particular, every element of $\ips{F}$ is the image of a uniformizer of $v_{t}$ under a polynomial over $F$.
\end{lemma}
\begin{proof}
If $a\in F$ then we may choose any $s\in\WW_{m}$ and choose $f$ to be the constant polynomial $a$.
Thus we may assume $a\notin F$.
Since $F=\bigcap_{k\geq0}F[\![t^{p^{k}}]\!]$,
there exists a unique $k\geq0$ such that
$a\in F[\![t^{p^{k}}]\!]\setminus F[\![t^{p^{k+1}}]\!]$.
First we suppose $k=0$,
so that
$a\notin F[\![t^{p}]\!]$,
and therefore
there exists $i\in\mathbb{N}$
such that
$a_{i}\neq0$ and $p\nmid i$.
We let
$w:=\min\{i\in\mathbb{N}:a_{i}\neq0\text{ and }p\nmid i\}\in\mathbb{N}$
and
$n:=m+2(w-1)$,
so that
$n\geq w$,
and we write
$f(X)=\sum_{i=0}^{n}a_{i}X^{i}$
and $g(X)=t^{m}(f(X)-a)$.
Note that
$v_{t}(g(t))>m+n$
and
$v_{t}(g'(t))=m+w-1$,
where $g'$ is the formal derivative of $g$.
Therefore
$v_{t}(g(t))>2v_{t}(g'(t))$.
By Hensel's Lemma, in the form of
{\cite[Theorem 4.1.3.(5)]{EP}},
since $g\in\mathcal{O}[X]$ and $t\in\mathcal{O}$,
there exists $s\in\mathcal{O}$ with
$g(s)=0$
and
$v_{t}(t-s)>v_{t}(g'(t))=m+w-1\geq m$.
In particular $s\in\WW_{m+1}\subseteq\WW_{m}$
satisfies $a=f(s)$.
Next suppose 
$k>0$.
Since $F$ is perfect, $a^{p^{-k}}\in\ips{F}\setminus F[\![t^{p}]\!]$,
and we may apply the case $k=0$ to the element $a^{p^{-k}}$ to find
$f=\sum_{i=0}^{n}b_{i}X^{i}\in F[X]$ and $s\in\WW_{m}$
such that $a^{p^{-k}}=f(s)$.
Applying the Frobenius endomorphism $k$ times, we find $a=f^{p^{k}}(s)$
where
$f^{p^{k}}$ denotes the polynomial $\sum_{i=0}^{n}b_{i}^{p^{k}}X^{ip^{k}}$.
The final claim follows from the inclusion $\WW_{m}\subseteq\UU$.
\end{proof}

\begin{remark}\label{rem:imperfect}
For any non-constant polynomial $f\in F[X]$ and any uniformizer $s\in\UU$,
the extension $F(s)/F(f(s))$ is algebraic and $\ps{F}/F(s)$ is separable
---
and when $F$ is perfect
every element of $\ips{F}\setminus F$ is of the form $f(s)$,
by Lemma~\ref{lem:Hensel_2}.
However,
if $F$ is imperfect there exists $a\in\ips{F}$ for which there do not exist $f\in F[X]$ and $s\in\UU$ such that $a=f(s)$.
To see this,
let $u\in F\setminus F^{(p)}$
and
let $b\in \ps{F}$ be transcendental over $F(t)$
---
such $b$ exist since $\ps{F}$ is transcendental over $F(t)$.
Then $a:=b^{p}+ut^{p}$
is an element of $\ps{F}$
such that there is no algebraic extension $E/F(a)$
such that 
$\ps{F}/E$ is separable,
since any such $E$ must contain both $b$ and $t$, so in particular has transcendence degree at least two over $F$.
Therefore 
the hypothesis that $F$ is perfect may not be removed from Lemma~\ref{lem:Hensel_2}.
\end{remark}

\subsection{$F$-automorphisms of $\ps{F}$}
\label{section:automorphisms}

Composition of formal power series (with respect to $t$) is the operation defined
as follows:
\begin{align*}
	\circ:\ps{F}\times(\mathfrak{m}\setminus\{0\})&\rightarrow\ps{F}\\
	(a,b)&\mapsto a\circ b:=\sum_{i\in\mathbb{Z}}a_{i}\big(\sum_{j>0}b_{j}t^{j}\big)^{i},
\end{align*}
where $a=\sum_{i\in\mathbb{Z}}a_{i}t^{i}$ and $b=\sum_{j>0}b_{j}t^{j}$.
It is easily verified that the above sum converges,
in the sense that the coefficient in $a\circ b$ of each $t^{i}$ is a finite sum,
so composition is well-defined by the above formula.
Each
$a\in\ps{F}$
induces a function
\begin{align*}
	\lambda_{a}:\mathfrak{m}\setminus\{0\}&\rightarrow\ps{F}\\
	b&\mapsto a\circ b,
\end{align*}
which we call {\em application} of $a$.
Similarly, each
$b\in\mathfrak{m}\setminus\{0\}$
induces a function
\begin{align*}
	\rho_{b}:\ps{F}&\rightarrow \ps{F}\\
	a&\mapsto a\circ b,
\end{align*}
which we call {\em substitution} by $b$.
Let $\mathbf{G}$ denote the group of $F$-automorphisms of $\ps{F}$,
i.e.~those field automorphisms that fix $F$ pointwise.
In \cite{S44}, Schilling gives a description of $\mathbf{G}$, summarized in the following fact.

\begin{fact}[{cf \cite[Theorem 1]{S44}}]\label{fact:Schilling}
Composition $\circ:\ps{F}\times(\mathfrak{m}\setminus\{0\})\rightarrow \ps{F}$ is continuous
with respect to the valuation topology.
The restriction of $\circ$ to $\UU\times\UU$ is associative, $t$ is the identity element, and every element of $\UU$ is invertible.
For each $b\in\mathfrak{m}\setminus\{0\}$, the substitution $\rho_{b}$ is a ring homomorphism, which fixes $F$ pointwise.
Thus $(\UU,\circ)$ is a group which acts on $\ps{F}$ as a group of $F$-automorphisms.
Moreover, the corresponding representation $\rho:(\UU,\circ)\rightarrow\mathbf{G}$, with $\rho(b)=\rho_{b}$, is an isomorphism.
\end{fact}

For $n\geq2$, let $\mathbf{G}_{n}$ denote the subgroup of $\mathbf{G}$ consisting of those automorphisms given by substitutions $\rho_{b}$, for $b\in\WW_{n}$.
The following chain of subgroups corresponds to
(\ref{eq:chain_1}):
\begin{align}\label{eq:chain_2}\tag{2}
		\mathbf{G}
	>
		\mathbf{G}_{2}
	>
		\mathbf{G}_{3}
	>
		\ldots
	>
		\bigcap_{n\geq2}\mathbf{G}_{n}
	=
		\{\mathrm{id}_{\ps{F}}\}.
\end{align}
For $n\geq2$, each $\WW_{n}$ is the orbit of $t$ under the action of $\GGG_{n}$.
Schilling proves in \cite[Theorem 3]{S44} that these groups are the same as the {pseudo-ramification groups} of MacLane, see \cite[Section 9]{ML39}.
We denote by $\Orb(\mathbf{a})$ the orbit of an $r$-tuple $\mathbf{a}=(a_{1},\ldots,a_{r})\in\ps{F}^{r}$ under $\GGG$,
and by $\Orb_{n}(\mathbf{a})$ the orbit of $\mathbf{a}$ under $\GGG_{n}$, for $n\geq2$.

\subsection{Suitable expansions}
\label{section:languages}

The first-order language of rings is $\Lring=\{+,\cdot,-,0,1\}$,
and naturally we may view each ring $R$ as an $\Lring$-structure, also denoted by $R$.
For any language $\mathfrak{L}\supseteq\Lring$ and any subset $A\subseteq R$,
$\mathfrak{L}(A)$ denotes the expansion of $\mathfrak{L}$ by an additional constant symbol $c_{a}$,
for each $a\in A$.
Then every $\mathfrak{L}$-structure of which $A$ is a subset may be expanded to an $\mathfrak{L}(A)$-structure in the natural way,
that is by interpreting $c_{a}$ as $a$, for each $a\in A$.

\begin{definition}\label{def:suitable}
An expansion $K=(\ps{F},\ldots)$ to an $\mathfrak{L}$-structure,
for a language $\mathfrak{L}\supseteq\Lring$,
is {\em suitable}
if the orbit of $t$ under the $\mathfrak{L}(F)$-automorphisms of $K$ contains $\mathcal{W}_{n}$, for some $n\geq2$.
\end{definition}

Straight from the definition, we observe that an $\mathfrak{L}$-expansion $K$ of $\ps{F}$ is suitable if and only if the natural further expansion to the language $\mathfrak{L}(F)$ is suitable.

There are many suitable expansions of $\ps{F}$, illustrated by the following example.

\begin{example}
\
\begin{enumerate}
\item
Of course $\ps{F}$ itself is suitable since $\UU$ is the orbit of $t$ under the group $\GGG$ of $\Lring(F)$-automorphisms of $\ps{F}$, and $\UU\supseteq\WW_{2}$.

\item
Let $\Lval=\Lring\cup\{O\}$ be the one-sorted language of valued fields,
where $O$ is a unary predicate symbol,
and let $K_{0}=(\ps{F},\mathcal{O}_{v_{t}})$ denote the expansion of $\ps{F}$ to an $\Lval$-structure obtained by interpreting $O$ by the $t$-adic valuation ring $\mathcal{O}=\ips{F}$.
Schilling shows in \cite[Lemma 1]{S44} that the $\Lring$-automorphisms of $\ps{F}$ are in fact already $\Lval$-automorphisms of $K_{0}$,
thus in particular $\GGG$ is equal to the group of $\Lval(F)$-automorphisms of $K_{0}$.
In particular, the expansion $K_{0}$ is suitable.
This can be seen rather easily from a model-theoretic point of view, since the valuation ring $\ips{F}$ is $\Lring$-definable in $\ps{F}$, without parameters, by a result of
Ax (\cite{Ax65}).
See also Robinson
(\cite{R65}).

\item
Let $\Lval^{3}$ be the three-sorted language of valued fields
with sorts
$\mathbb{K}$, $\mathbbm{k}$, and $\mathbb{G}$
equipped with $\Lring$, $\Lring$, and the language $\mathfrak{L}_{\mathrm{oag}}\cup\{\infty\}$ of ordered abelian groups with extra symbol $\infty$, respectively,
as well as two function symbols $v:\mathbb{K}\rightarrow\mathbb{G}$ and $\res:\mathbb{K}\rightarrow\mathbbm{k}$ for the valuation and residue maps, respectively.
The expansion $K_{1}$ of $\ps{F}$ to an $\Lval^{3}$-structure may be defined in the obvious way,
and by setting the residue map to $0$ outside of the valuation ring.
The restriction of each $\Lval^{3}$-automorphism of $K_{1}$ to its action on the sort $\mathbb{K}$
defines a bijection between the $\Lval^{3}$-automorphisms of $K_{1}$ and the $\Lval$-automorphisms of $K_{0}$.
In particular, the expansion $K_{1}$ is suitable.

\item
Let $\Lac$ denote 
the expansion of $\Lval$ by a unary function symbol $\mathrm{ac}$.
Let $K_{2}$ denote the expansion of $K_{0}$ by the interpreting $\mathrm{ac}$ by the {\em angular component}, which is the map
\begin{align*}
	\underline{\mathrm{ac}}:\ps{F}^{\times}&\rightarrow F^{\times}\\
	a=\sum_{i\geq n}a_{i}t^{i}&\mapsto a_{v_{t}(a)}
\end{align*}
The group of $\Lac(F)$-automorphisms of $K_{2}$ is
is equal to $\mathbf{G}_{2}$,
thus $K_{2}$ is suitable.

\item
The {\em simple analytic functions} on $\ps{F}$ are the maps
$\lambda_{a}$, for $a\in\ips{F}$, extended by mapping $\lambda_{a}(0):=a_{0}$.
Let $\Lan$ denote 
the language 
$\Lring$ expanded by unary function symbols $f_{a}$,
and let $K_{3}$ denote the expansion of $\ps{F}$ to an $\Lan$-structure
by interpreting each $f_{a}$ by the corresponding simple analytic function $\lambda_{a}$, extended as described above.
In fact $K_{3}$ has the same group of automorphisms as $\ps{F}$, namely $\GGG$,
thus it is suitable.

\item
For each uniformizer $s\in\UU$,
there is a {\em cross-section}
$\chi_{s}:\mathbb{Z}\rightarrow\ps{F}$,
$n\mapsto s^{n}$,
viewed as a map with domain the value group.
The family
$(\chi_{s})_{s\in\UU}$
of cross-sections
may be combined into one binary function 
$\chi:\mathbb{Z}\times\ps{F}\rightarrow\ps{F}$
that maps $(n,s)$ to $\chi_{s}(n)$,
in case $s\in\UU$,
and takes (for example) the value $0$ otherwise.
The expansion $(K_{1},\chi)$ of $K_{1}$ by $\chi$
is one way of formalizing the expansion of $\ps{F}$ by the parameterized family $(\chi_{s})_{s\in\UU}$.
This structure has the same group of automorphisms as $\ps{F}$,
therefore it is suitable.
\end{enumerate}
\end{example}

On the other hand, several expansions of $\ps{F}$ are rather rigid, having groups of $F$-automorphisms not containing $\GGG_{n}$, for any $n\geq2$.

\begin{example}
\
\begin{enumerate}
\item
There are no non-trivial $\mathfrak{L}_{\mathrm{ring}}(F(t))$-automorphisms of $\ps{F}$.
Thus $(\ps{F},t)$ is not suitable.
\item
Consider the single cross-section $\chi_{t}:\mathbb{Z}\rightarrow\ps{F}$ that maps $n\mapsto t^{n}$.
The structure $(\ps{F},\chi_{t})$ has no non-trivial $F$-automorphisms,
thus it is not suitable.
\end{enumerate}
\end{example}

\subsection{Definability}

When stating a result on first-order definability, we will be explicit about which language and set of parameters are allowed.
That is,
we will write e.g.~``{\em $\mathfrak{L}(A)$-definable}'' to mean that the definition is with an $\mathfrak{L}$-formula and parameters drawn from $A$, which would be a subset of the domain of the structure under consideration.
Otherwise we will write ``{\em parameter-definable}'' to mean that we may use any parameters from the given model.

An $\mathfrak{L}$-formula is {\em existential} if it is of the form $\exists y_{1}\ldots\exists y_{n}\;\psi(x_{1},\ldots,x_{m},y_{1},\ldots,y_{n})$ for a quantifier-free $\mathfrak{L}$-formula $\psi$, and sets defined by existential formulas are {\em existentially} definable.

\begin{lemma}\label{lem:definitions}
The set\;$\UU$ is
(a) $\Lring$-definable and
(b) existentially $\Lring(t)$-definable.
The sets $\WW_{n}$, for $n\geq2$, are
existentially $\Lring(t)$-definable.
\end{lemma}
\begin{proof}
The valuation ring $\mathcal{O}=\ips{F}$ is definable in $\ps{F}$
by an existential $\Lring(t)$-formula $\psi(x,t)$,
due to Robinson
(\cite{R65}).
For convenience, we write $\psi(x,t)$ with the parameter $t$ made explicit.
Since $\UU=t(\mathcal{O}\cap(\mathcal{O}\setminus\{0\})^{-1})$, the set of uniformizers is defined by the
formula
$$
	\exists y\exists z\;(\psi(y,t)\wedge\psi(z,t)\wedge yz=1\wedge x=ty),
$$
which is logically equivalent to an existential $\Lring(t)$-formula, proving (b).

Next we recall that there is an $\Lring$-formula $\psi'(x)$ due to Ax (\cite{Ax65}) that defines $\ips{F}$ in $\ps{F}$.
The set of units $\mathcal{O}^{\times}$ is defined by the formula $\chi'(x)$ given by
$\psi'(x)\wedge\exists y\;(xy=1\wedge\psi'(y))$,
and the maximal ideal $\mathfrak{m}$ is defined by the formula $\xi'(x)$ given by
$\psi'(x)\wedge\neg\chi'(x)$.
Finally, $\UU$ is defined by the formula
\begin{align*}
	\psi'(x)\wedge\neg\chi'(x)\wedge\neg\exists y\exists z(x=yz\wedge\xi'(y)\wedge\xi'(z)).
\end{align*}
The final claim is clear from (b) and the description of $\WW_{n}$ as $t+t^{n}\ips{F}$.
\end{proof}

\begin{proposition}\label{prp:main}
Suppose $F$ is perfect.
Let $K=(\ps{F},\ldots)$ be a suitable $\mathfrak{L}$-expansion of $\ps{F}$,
and let $a\in K\setminus F$.
Then the orbit of $a$ under the $\mathfrak{L}$-automorphisms of $K$ contains an
infinite existentially $\Lring(F(t))$-definable set.
In particular,
if $X\subseteq \ps{F}$
is $\mathfrak{L}(F)$-definable
and $X\not\subseteq F$,
then there is an infinite subset $Y\subseteq X\setminus F$ that is existentially $\Lring(F(t))$-definable in $\ps{F}$.
\end{proposition}
\begin{proof}
Since $K$ is suitable, there exists $n\geq2$ such that the orbit of $t$ under the $\mathfrak{L}(F)$-automorphisms of $K$ contains $\WW_{n}$.
If $a\in\ips{F}$,
then
$a=f(s)$ for some $f\in F[X]$ and $s\in\WW_{n}$,
by Lemma~\ref{lem:Hensel_2}.
In particular, $s$ and $t$ lie in the same orbit under $\mathfrak{L}(F)$-automorphisms.
The orbit of $a$ under the $\mathfrak{L}(F)$-automorphisms therefore contains $f(\WW_{n})$,
which by Lemma~\ref{lem:definitions}
is existentially $\Lring(F(t))$-definable.
Such a set is clearly infinite since it contains a transcendence basis of $\ps{F}/F$.
On the other hand, if $a\notin\ips{F}$ then $a^{-1}\in t\ips{F}$.
By what we have just shown, the orbit of $a^{-1}$ under $\mathfrak{L}(F)$-automorphisms
contains the infinite set defined by some existential $\Lring(F(t))$-formula
$\varphi(x)$.
The orbit of $a$ under the same group of automorphisms thus contains the set defined by
$\exists y\;(xy=1\wedge\varphi(y))$,
which is logically equivalent to an existential $\Lring(F(t))$-formula.
The final assertion follows from considering the orbit under $\mathfrak{L}(F)$-automorphisms of any $a\in X\setminus F$.
\end{proof}

We take this opportunity to recall some notation from \cite[\S3]{A19}.

\begin{definition}
Fix an enumeration $(f_{i})_{i<\omega}$ of the multivariable polynomials over $\mathbb{Z}$.
We may arrange the enumeration so that each $f_{i}$ is a polynomial in (at most) the variables $X_{0},...,X_{i-1}$.
For any formula $\varphi(x)$ in one-free variable and any $m<\omega$,
we let $\varphi_{m}(x)$ denote the formula
\begin{align*}
\exists\;
\mathbf{a}=(a_{k})_{k<m},\mathbf{b}=(b_{l})_{l<m}\;:\;
\bigvee_{i,j<m}\bigg(
x\cdot f_{j}(\mathbf{a})=f_{i}(\mathbf{b})
\wedge
f_{j}(\mathbf{a})\neq0
\wedge\bigwedge_{k,l<m}\Big(\varphi(a_{k})\wedge\varphi(b_{l})\Big)\bigg).
\end{align*}
\end{definition}

Note that if $\varphi(x)$ is an existential $\mathfrak{L}$-formula,
then $\varphi_{m}(x)$ is also (logically equivalent to) an existential $\mathfrak{L}$-formula, for any language $\mathfrak{L}\supseteq\Lring$.
In an $\mathfrak{L}$-expansion $K$ of a field, $\varphi_{m}(x)$ defines the increasing union of images of the set defined by $\varphi(x)$ under the rational functions $\frac{f_{i}}{f_{j}}$, for $i,j<m$.
For a set $X$ defined by a formula $\varphi(x)$, we denote by $X_{m}$ the set defined by $\varphi_{m}(x)$.
Thus the field $(X)$ is the increasing union $\bigcup_{m<\omega}X_{m}$.

\begin{theorem}\label{thm:subfields}
Suppose $F$ is perfect.
Let $K=(\ps{F},\ldots)$ be a suitable $\mathfrak{L}$-expansion of $\ps{F}$,
and let $X\subseteq \ps{F}$ be an $\mathfrak{L}(F)$-definable subset.
Then either $X\subseteq F$ or there exists $m,n\geq0$ such that
$(X)=X_{m}=F(\!(t^{p^{n}})\!)$.
\end{theorem}
\begin{proof}
If $X\not\subseteq F$ then by Proposition~\ref{prp:main}
there is an infinite existentially $\Lring(F(t))$-definable set $Y\subseteq X$.
By~{\cite[Proposition 30]{A19}},
there exist $m,n\geq0$ such that
$(Y)_{m}\supseteq F(\!(t^{p^{n}})\!)$.
Thus $(X)_{m}$ contains $F(\!(t^{p^{n}})\!)$.
The conclusion follows from the fact that
the subfields intermediate between $F(\!(t^{p^{n}})\!)$ and $\ps{F}$ are exactly those fields $F(\!(t^{p^{k}})\!)$ for $k$ such that $0\leq k\leq n$.
\end{proof}

Theorem~\ref{thm:subfields_intro} is an immediate corollary of Theorem~\ref{thm:subfields} in the case $\mathfrak{L}=\Lval(F)$.
In fact, in the proof of Theorem~\ref{thm:subfields_intro} we may appeal to
{\cite[Theorem 1]{A19}} instead of \cite[Proposition 30]{A19} since we have made no claim about $X_{m}$.

There are three rather simple further corollaries of this theorem.

\begin{corollary}\label{cor:dichotomy}
Suppose $F$ is perfect
and let $X$ be an $\Lval(F)$-definable subset of $\ps{F}$.
Then either $X\subseteq F$ or $F\subseteq(X)$.
\end{corollary}

\begin{corollary}\label{cor:finite.constants}
Suppose $p$ is a prime number and let $q=p^{k}$ be a prime power,
let $K=(\ps{\mathbb{F}_{q}},\ldots)$ be a suitable $\mathfrak{L}$-expansion of $\ps{\mathbb{F}_{q}}$,
and
let $X$ be an $\mathfrak{L}(F)$-definable subset of $K$.
Then $(X)$ is existentially $\Lring$-definable.
\end{corollary}
\begin{proof}
It follows from Theorem~\ref{thm:subfields} that $(X)$ is
either a subfield of $\mathbb{F}_{q}$
or of the form
$F(\!(t^{p^{n}})\!)$, for some $n<\omega$.
In either case, $(X)$ is existentially $\Lring$-definable.
\end{proof}

\begin{corollary}\label{cor:perfect_hull}
Suppose $F$ is perfect and
let $X$ be an $\Lval(F)$-definable subset of $\ps{F}^{\mathrm{perf}}$.
Then either $X\subseteq F$ or $(X)=\ps{F}^{\mathrm{perf}}$.
\end{corollary}
\begin{proof}
Suppose $X\not\subseteq F$ and let $a\in X\setminus F$.
Let $n\in\mathbb{Z}$ be chosen maximal such that
$a\in F(\!(t^{p^{n}})\!)$,
so that
$a\notin F(\!(t^{p^{n+1}})\!)$.
We write $s:=t^{p^{n}}$
and observe that each
$\Lval(F)$-automorphism of $F(\!(s)\!)$ extends (uniquely) to an $\Lval(F)$-automorphism of $F(\!(s)\!)^{\mathrm{perf}}=\ps{F}^{\mathrm{perf}}$.
Thus the set
$X_{s}:=X\cap F(\!(s)\!)$ is invariant under
$\Lring(F)$-automorphisms of $F(\!(s)\!)$,
and in particular it contains the orbit of $a$ under these automorphisms.
In turn,
by Proposition~\ref{prp:main},
this orbit contains
an infinite existentially $\Lring(F(s))$-definable subset $Y$ of $F(\!(s)\!)$.
We again apply \cite[Proposition 30]{A19}
to find $m,n<\omega$
such that
$Y_{m}\supseteq F(\!(s^{p^{n}})\!)$.
Therefore
$(X_{s})\supseteq X_{m}\supseteq Y_{m}\supseteq F(\!(s^{p^{n}})\!)$;
and since $a\in Y$ is a $p$-basis of $F(\!(s)\!)$,
we even have
$(X_{s})=F(\!(s)\!)$.
Finally we consider the automorphism $f$ of $\ps{F}^{\mathrm{perf}}$ that fixes $F$ pointwise and sends $t\mapsto t^{1/p}$.
The set $X$, and thus the field $(X)$, are both closed under $f$,
which yields
$(X)=\ps{F}^{\mathrm{perf}}$, as required.
\end{proof}

\begin{remark}
Since 
$\ps{F}^{\mathrm{perf}}$ is a perfect large field,
it seems likely that one can give an alternative proof of Corollary~\ref{cor:perfect_hull} in which \cite[Proposition 30]{A19} is replaced by a result from \cite{Fehm}.
\end{remark}

\begin{remark}
These results can be seen in the context of \cite[Corollary 5.6]{JK10}, in which it is shown that a henselian field of characteristic zero has no proper parameter-definable subfields (in the language of rings);
and \cite[Question 10]{JK10} in which Junker and Koenigsmann ask whether $\ps{\mathbb{F}_{p}}^{\mathrm{perf}}$ is very slim (see {\cite[Definition 1.1]{JK10}}).
If $\ps{\mathbb{F}_{p}}^{\mathrm{perf}}$ were very slim then in particular it would have no infinite proper parameter-definable subfields.
Corollary~\ref{cor:perfect_hull} shows that $\ps{\mathbb{F}_{p}}^{\mathrm{perf}}$ has no infinite proper subfields which are definable with parameters only from $\mathbb{F}_{p}$,
but at present we are unable to extend our methods to study subfields definable with parameters from outside the subfield of constants.
\end{remark}

\section{Orbits of one-dimensional tuples}
\label{section:1-dim_orbits}

The {henselian valuation topology} on $\ps{F}$ is the non-discrete Hausdorff topology induced by $v_{t}$:
by definition, the sets 
$\B(n;a)=\{x\in\ps{F}\mid v_{t}(x-a)>n\}$,
for $n\in\mathbb{Z}$,
form a basis for the filter of open neighbourhoods of each $a\in\ps{F}$.
Note that $\B(n;t)=\WW_{n+1}$ for every $n<\omega$.

\begin{theorem}\label{thm:1-dim}
\perfect\
Let $\mathbf{a}$ be a tuple from $\ps{F}$ of transcendence degree $1$ over $F$.
Then $\Orb(\mathbf{a})$ is
\begin{enumerate}
\item existentially $\Lring(F(t))$-definable,
\item $\Lring(F)$-definable, and 
\item the set of realisations in $\ps{F}$ of the $\Lring$-type $\mathrm{tp}(\mathbf{a}/F)$.
\end{enumerate}
\end{theorem}
\begin{proof}
By replacing each element of $\mathbf{a}$ by its multiplicative inverse, if necessary, we may assume that $\mathbf{a}$ is a tuple from $\ips{F}$.
For example,
to see that the truth of \I\ does not alter under such a replacement:
if $\mathbf{a}$ is partitioned as $(b,\mathbf{c})$, with $b\neq0$,
and if an existential $\Lring(F(t))$-formula $\varphi(\mathbf{x})$ defines $\Orb(\mathbf{a})$ then $\exists y'(\varphi(y',\mathbf{z})\wedge yy'=1)$ is logically equivalent to an existential $\Lring(F(t))$-formula that defines $\Orb(b^{-1},\mathbf{c})$.

By reordering if necessary, we may write $\mathbf{a}$ as $(b,c_{1},\ldots,c_{r})$ such that $b$ is transcendental over $F$ and each $c_{i}$ is algebraic over $F(b)$.
By Lemma~\ref{lem:Hensel_2}, there is a polynomial $f\in F[X]$ and a uniformizer $s\in\UU$ such that $f(s)=b$.
Because the henselian valuation topology is Hausdorff,
there exists $n\in\mathbb{N}$ such that each $c_{i}$ in $\mathbf{c}$ is the unique zero in $\B(n;c_{i})$ of its minimal polynomial over $F(b)$.
It follows that there is a quantifier-free $\Lring(F(s))$-formula $\varphi_{i}(z,s)$ such that the intersection of the set it defines in $\ps{F}$ and $\B(n;c_{i})$ is $\{c_{i}\}$.
Next, the ring $F[t]$ is dense in $\ips{F}$,
with respect to the henselian valuation topology,
and so $F[s]$ is also dense in $\ips{F}$ by symmetry.
Therefore there is a polynomial $g_{i}\in F[X]$ such that $g_{i}(s)\in\B(n;c_{i})$, and by applying the ultrametric inequality we have $\B(n;c_{i})=g_{i}(s)+s^{n+1}F[\![s]\!]$.
By Robinson (\cite{R65}), there is an existential $\Lring(s)$-formula $\psi(w,s)$ that defines $F[\![s]\!]$ in $\ps{F}$.
Thus the ball $\B(n;c_{i})$ itself is $\Lring(s)$-definable
by the formula $\chi_{i}(z,s)$, defined to be
$\exists w\;(\psi(w,s)\wedge z=g_{i}(s)+s^{n+1}w)$,
which is logically equivalent to an existential $\Lring(F(s))$-formula.
Therefore
\begin{align*}
\Big(f(s)=y\wedge\bigwedge_{i}(\chi_{i}(z_{i},s)\wedge\varphi_{i}(z_{i},s))\Big)
\end{align*}
is an $\Lring(F(s))$-formula that defines the single tuple $\mathbf{a}=(b,\mathbf{c})$.
Combining this formula with 
the existential $\Lring(s)$-definition (respectively, the $\Lring$-definition)
of $\UU$,
as in Lemma~\ref{lem:definitions},
we obtain an existential $\Lring(s)$-definition (respectively, an $\Lring$-definition)
of the orbit of $\mathbf{a}$ under the action of $\GGG$.
This proves \I\ and \II.
For \III\ we first note that $\Orb(\mathbf{a})$ is a subset of the set of realisations in $\ps{F}$ of the $\Lring(F)$-type of $\mathbf{a}$, which itself is the intersection of all $\Lring(F)$-definable sets containing $\mathbf{a}$.
Since even $\Orb(\mathbf{a})$ is already $\Lring(F)$-definable, by \II, this proves that $\Orb(\mathbf{a})$ coincides with the set of realisations of the $\Lring(F)$-type of $\mathbf{a}$.
\end{proof}

Theorem~\ref{thm:1-dim_intro}
is immediate from Theorem~\ref{thm:1-dim}.

\section*{Acknowledgements}

This research is based on a part of the author's doctoral thesis
completed under the supervision of Jochen Koenigsmann and supported by EPSRC.
The author would like to thank
Arno Fehm and Franziska Jahnke
for invaluable feedback on earlier versions.

%%%Bibliography%%%
\def\bibfont{\footnotesize}
\bibliographystyle{plain}
% \bibliography{bibliography}
\begin{singlespacing}

\end{singlespacing}
\end{document}